\documentstyle[12pt]{article}
\title{Noncommutative analogues of q-special polynomials and q-integral 
on a quantum sphere}
\author{D. Gurevich* and L. Vainerman**}
\date{}
\begin{document}
\newtheorem{theorem}{Theorem}
\newtheorem{lemma}{Lemma}
\newtheorem{claim}{Claim}
\newtheorem{corollary}{Corollary}
\newtheorem{proposition}{Proposition}
\newtheorem{definition}{Definition}
\newtheorem{conjecture}{Conjecture}
\newtheorem{remark}{Remark}
\newcommand{\ren}{\rho^{End}_q}
\newcommand{\orr}{\overline{\rho}}
\newcommand{\br}{[\,\,,\,\,]}
\newcommand{\brq}{[\,\,,\,\,]_q}
\newcommand{\vqp}{V^q_+}
\newcommand{\vqm}{V^q_-}
\newcommand{\g}{\bf g}
\newcommand{\ahq}{A^c_{h,q}}
\newcommand{\aaa}{A^c_{0,1}}
\newcommand{\aq}{A^c_{0,q}}
\newcommand{\ah}{A^c_{h,1}}
\newcommand{\la}{\lambda}
\newcommand{\Int}{{\rm Int}}
\newcommand{\tS}{\widetilde{S}}
\newcommand{\ahc}{A_{h,1}^c}
\newcommand{\rn}{\rho_{\nu}}
\newcommand{\usl}{U_q(sl(2))}
\newcommand{\gggg}{g}
\newcommand{\ug}{U(\gggg)}
\newcommand{\us}{U(sl(2))}
\newcommand{\vv}{{\bf V}^{\otimes 2}}
\newcommand{\ww}{W^{\otimes 2}}
\newcommand{\uqs}{U_q(sl(2))}
\newcommand{\ork}{\overline{\rho_k}}
\newcommand{\uog}{U(\overline{\gggg})}
\newcommand{\uos}{U(\overline{sl(2)})}
\newcommand{\uq}{U_q(\gggg)}
\newcommand{\sn}{sl(n)}
\newcommand{\ogg}{\overline{\gggg}}
\newcommand{\osl}{\overline{sl(2)}}
\newcommand{\osln}{\overline{sl(n)}}
\newcommand{\rd}{\rho^{\otimes 2}}
\newcommand{\vbq}{V_{\beta}^q}
\newcommand{\uqq}{U^q_{1/2}}
\newcommand{\uqk}{U^q_{k/2}}
\newcommand{\vq}{V^q}
\newcommand{\vb}{V_{\beta}}
\newcommand{\vqb}{V_{\beta}^q}
\newcommand{\voq}{V_{\omega}^q}
\newcommand{\vkoq}{V_{k\omega}^q}
\newcommand{\vkb}{V_{k\beta}}
\newcommand{\vbk}{V_{\beta}^{\otimes k}}
\newcommand{\vko}{V_{k\omega}}
\newcommand{\vlo}{V_{l\omega}}
\newcommand{\vklo}{V_{(k+l)\omega}}
\def\ot{\otimes}
\def\De{\Delta}
\def\qq{q^{-1}}
\maketitle
\centerline {* ISTV, Universit\'e de Valenciennes, 59304 Valenciennes, France,}
\centerline{e.mail: gurevich@univ-valenciennes.fr}
\centerline {** International Solomon University, Zabolotny Street, 38, apt. 61,}
\centerline {Kiev 252187 Ukraine, e.mail: wain@agrosys.kiev.ua}
\begin{abstract}

The q-Legendre polynomials can be treated as some special "functions in the quantum double 
cosets $U(1)\setminus SU_q(2)/U(1)$". They form a family (depending on a parameter $q$) 
of polynomials in one variable. We get their further generalization by introducing 
a two parameter family of polynomials. If the former family arises from an algebra 
which is in a sense "q-commutative", the latter one is related to its 
noncommutative counterpart. We introduce also a two parameter deformation of 
the invariant integral on a sphere. 

\end{abstract}

\section{Introduction}

It is well known that the classical Legendre polynomials form a basis in the 
function space on the double cosets $U(1)\setminus SU(2)/U(1)$ [Vi]. Although 
q-analogues of these polynomials (as well as those of some other special 
functions) are known for a long time, it became clear only recently  that they
can be treated as "functions on quantum double cosets".

This approach suggested in \cite{VS} for the $sl(2)$ case (see also \cite{KN}) and 
developed by a number of authors for other quantum double cosets (cf. the survey 
\cite{Va}) can be explained as follows.  Let us consider the function space $Fun_q(S^2)$ 
on the quantum sphere. This space can be defined in the spirit of the paper \cite{P} as 
the subspace of left (or right) $U(1)$-invariant functions on $SU_q(2)$ (the group $U(1)$
can be treated as a commutative and cocommutative Hopf subalgebra of $SU_q(2)$).

For a generic $q$ the space $Fun_q(S^2)$ can be decomposed into  a direct sum 
$\oplus {\bf V}_i,\,i=0,1,...$ of irreducible $U_q(su(2))$-modules ${\bf V}_i$, where $i$ 
is the spin (note that ${\rm dim} {\bf V}_i=2i+1$). Let $v$ be a generator of the subalgebra of 
$Fun_q(S^2)$ formed by two-sided $U(1)$-invariant functions. Then the k-th q-Legendre 
polynomial can be defined as a polynomial in $v$ belonging to the component ${\bf V}_k$. 
In fact, the q-Legendre polynomials are nothing but eigenfunctions of the quantum Casimir 
operator (this defines them uniquely up to factors). 

It should be noted that the algebra $Fun_q(S^2)$, which plays a crucial role in the 
constructions of \cite{VS}, is a particular case of a two parameter family of 
$U_q(su(2))$-invariant (in the sense explained in Section 2) associative algebras (meanwhile, 
the algebra $Fun_q(S^2)$ itself depends only on the  parameter $q$ assuming that the 
parameter $c$ labeling the orbits is fixed, cf. below). 

This two parameter family arises from a quantization of the Poisson pencil generated by 
the Kirillov-Kostant-Souriau (KKS) bracket on the usual sphere and the so-called R-matrix 
bracket (cf. the last Section for the definition).

The quantization of the  R-matrix bracket leads to the algebra $Fun_q(S^2)$ which plays 
the role of "commutative algebra" in the category of $\usl$-invariant algebras. The passage 
to the  two parameter family mentioned above is a  way to a "q-noncommutative" analysis. 
Our main aim is to apply the above approach to this family. More precisely, we introduce  
$(\hbar, q)$-special polynomials as eigenfunctions of the quantum Casimir operator acting 
on this two parameter family.

Besides,  we introduce in the spirit of the paper \cite{NM}, a certain $(\hbar, q)$-analogue of 
the invariant integral on the sphere (the authors of \cite{NM} deal with a $q$-analogue of the 
integral wich is defined on the algebra $Fun_q(S^2)$). And finally, we give an explicit expression 
for this $(\hbar,q)$-integral which is a 
generalization of the well known Jackson integral (see, for example \cite{VS}, \cite{KN}).

Throughout the paper, the basic field is $\bf C$ and $q\in\bf C$ is assumed to be generic.
Thus, we deal with the group $SL(2,{\bf C})$, the complexification of the sphere and their 
quantum counterparts  rather then with compact objects themselves: the reason for this is 
explained in Section 4.

The paper is organized as follows. In the following Section we define our basic object: a 
function algebra on a quantum hyperboloid. In Section 3 we compute the action of the 
quantum Casimir operator on the elements of this algebra.

In Section 4 we define  $(\hbar, q)$-special polynomials in the above algebra.

Section 5 is devoted to introducing an $(\hbar, q)$-analogue of the invariant integral on 
a sphere. In the last Section we discuss the considered objects in terms of deformation 
quantization. There we also explain in what sense we use the term {\em a q-commutative algebra}. 

\section{Basic objects: the algebra $\ahq$}

Let us consider the quantum enveloping algebra $\uqs$, i.e., an algebra generated by elements 
$ E_{+},\ E_{-},\ X,\ Y$ satisfying the following relations: 
$$E_{\pm}X=q^{\pm 1}XE_{\pm};\  E_{\pm}Y=q^{\mp 1}Y E_{\pm};\  
E_{+} E_{-}=E_{-}E_{+}=1;\ [X,Y]=\frac{ E_{+}^2 - E_{-}^2}{q-q^{-1}},$$
where $q\neq 0,q^2\neq 1,$ equipped with the coproduct
$$\De(X)= E_{-} \ot X+X\ot  E_{+} ;\ \De(Y)=E_{-} \ot Y+Y\ot  E_{+} ;
\De(E_{\pm})= E_{\pm}\ot E_{\pm}$$
and some antipode whose explicit form we do not need.

One can verify that the element, called {\it quantum Casimir operator} or simply
quantum Casimir,
$$K=\frac{q}{2}(XY+YX)+\frac{q^2(1+q^2)}{2(1-q^2)^2}(E_{+}^2+ E_{-}^2 -2)$$
belongs to the center of the above algebra.

It is well known that the image of the quantum Casimir in any irredicible 
$\uqs$-module is a scalar operator. Let us denote by $\la_k,\, k=0,1/2,1, ...$
its eigenvalue corresponding to an  irreducible spin $k$ $\uqs$-module
${\bf V}_k$. We will show later that 
\begin{equation}
\la_k=\frac{(q^{2k}-1)(q^{2(k+1)}-1)}{q^{2k-2}(q^2-1)^2}.
\end{equation}

Now, let us consider a three-dimensional $\uqs$-module ${\bf V}={\bf V_1}$ 
such that the representation $\rho_q: \uqs\to End({\bf V})$ coincides with the
classical one  $\rho: U(sl(2))\to End({\bf V})$  as $q=1$. Let us fix a basis 
$\{u,\,v,\,w\}$ in  ${\bf V}$ such that the above action of the quantum group 
is given by (we omit the symbol $\rho_q$ in our notations): 
$$E_{\pm}u=q^{\pm 1}u,\ E_{\pm}v=v,\ E_{\pm}w=q^{\mp 1}w,\ Xu=0,\ Xv=-(q+\qq) u,\
Xw=v,$$ $$Yu=-v,\ Yv=(q+\qq) w,\ Yw=0.$$

Using the coproduct, we can equip $\vv$ with a $\uqs$-module structure as well.
This module is reducible and can be decomposed into a direct sum of three
irreducible $\uqs$-modules
$${\bf V}_0 =span\{(q^2+1)uw+vv+\frac{q^2+1}{q^2} wu\},$$
$${\bf V}_1 =span\{q^2 uv - vu,\ (q^2+1)( uw-wu) +(1-q^2)vv,\ -q^2 vw + wv\},$$
$${\bf V}_2 =span\{uu,\ uv+q^ 2vu,\ \qq uw-qvv+q^3 wu,\ vw+q^2 wv,\ ww\}$$
of spins 0,1,2 respectively (here the sign $\ot$ is omitted).

Then only the following relations imposed on the elements of the space
$\vv\oplus {\bf V}\oplus {\bf C}$ are consistent with the above action of
$\uqs$: $$C_q =(q^2+1)uw+vv+\frac{q^2+1}{q^2} wu =c,$$
$$q^2 uv - vu=-\hbar u,$$
$$(q^2+1)( uw-wu) +(1-q^2)vv=\hbar v,$$
$$-q^2 vw + wv=\hbar w$$
with arbitrary $\hbar$ and $c$. The element $C_q$ is called {\it a braided Casimir}.

Let us denote $\ahq$ the quotient algebra of a free tensor algebra $T({\bf V})$ by 
the ideal generated by the elements 
$$(q^2+1)uw+vv+\frac{q^2+1}{q^2} wu -c,\ q^2 uv - vu+\hbar u,$$
$$(q^2+1)( uw-wu) +(1-q^2)vv-\hbar v,\ -q^2 vw + wv-\hbar w.$$

\begin{remark} Here $\hbar$ and $q$ are assumed to be fixed. If we want to consider them as 
formal parameters, we must replace $T({\bf V})$ in the definition of the algebra $\ahq$ by 
$T({\bf V})\ot{\bf C}[[\hbar, q]]$. The parameter $c$ which labels the orbits is always fixed. 
The case $c=0$ corresponds to the cone.
\end{remark}

If $q=1, \,\hbar =0$, we get a family (parametrized by the parameter $c$ which 
labels the orbits) of usual hyperboloids considered as orbits in $sl(2)^*$ (the case $c=0$
corresponds to the cone). If 
$q=1, \,\hbar \not=0$, we get its noncommutative analogue but it still lives 
in the classical category of $sl(2)$-invariant algebras.

If $q\not=1$, we get a two parameter family of $\usl$-invariant algebras. Let us recall 
that an associative algebra $A$ is called $U_q(\gggg)$-invariant (or covariant) if
$$X\circ(a\ot b)=\circ\,\De\, X\,(a\ot b),\,\,\forall X\in U_q(\gggg),\,a,b\in A$$
where $\circ$ is the product in $A$.

In fact, the Podles' quantum spheres are exactly these quantum hyperboloids equipped with an 
involution. Here, we would like to avoid a discussion of the problem of a proper definition 
of an involution in braided categories (it has been discussed in \cite{DGR1}) and prefer 
to work with complex objects.

The particular case $\hbar=0$ of this family corresponds to a q-commutative algebra in the
sense discussed in the last Section.

Let us rewrite the above equations as follows: 
$$(q^2+1)uw+{\tilde v}^2+\frac{q^2+1}{q^2} wu =\tilde c -2a\tilde v,$$
$$q^2 u\tilde v - \tilde v u=0,$$
$$(q^2+1)( uw-wu) +(1-q^2){\tilde v}^2=-\hbar\tilde v,$$
$$-q^2 \tilde v w + w\tilde v =0,$$
where $a=\hbar(1-q^2)^{-1},\ \tilde c =c-a^2,\tilde v =v-a$.

Using these relations we can express the product $uw$ in terms of the variable 
$\tilde v$:
\begin{equation}
uw=(q^2+1)^{-2}[{\tilde c}q^2 -a(1+q^2){\tilde v}-{\tilde v}^2].
\end{equation}

We will use this formula below.

\section{Action of the quantum Casimir}

Our next aim is to get a formula for $K{\tilde v}^k$ for every natural $k$,
where ${\tilde v}^k={\tilde v}^{\ot k}$.

It is clear from the very beginning that the action of
$E_{+}^2+ E_{-}^2 -2$ on ${\tilde v}^k$ equals to
$0$. On the other hand, from the relation of commutation for $X,Y$ it is also
clear that the actions of $XY$ and $YX$ on ${\tilde v}^k$ coincide. 

Thus, we have 
$$K {\tilde v}^k=q YX {\tilde v}^k.$$

Using the formulae for the coproduct and for the action of $X,E_{+},E_{-}$ on
${\tilde v}$ as well as the formula of commutation of $v$ and ${\tilde v}$, we
have: 
$$X{\tilde v}^k=(XE_{+}^{k-1}+E_{-}XE_{+}^{k-2}+...+E_{-}^{k-1}X){\tilde v}^k=$$
$$=-\frac{q^2+1}{q}(u{\tilde v}^{k-1} +{\tilde v}u{\tilde v}^{k-2}+...+
{\tilde v}^{k-1}u)=$$
$$=-\frac{q^2+1}{q}(1+q^2 +...+q^{2(k-1)})u{\tilde v}^{k-1}=
-\alpha_k (q)u{\tilde v}^{k-1}$$
with 
$$\alpha_k(q)=\frac{(q^2 +1)(q^{2k}-1)}{q(q^2 -1)}.$$

Hereafter by $XE_{+}^{k-1}{\tilde v}^k$ we mean $X{\tilde v}(E_{+}{\tilde v})^{k-1}$ etc.

Similarly, using the formula (2) we get:
$$YX{\tilde v}^k=-\alpha_k (q)(YE_{+}^{k-1}+E_{-}YE_{+}^{k-2}+...+E_{-}^{k-1}Y)
u{\tilde v}^{k-1}=$$
$$=-\alpha_k (q)[-({\tilde v}+a){\tilde v}^{k-1}+\frac{q^2+1}{q^2}
(uw{\tilde v}^{k-2}+u{\tilde v}w{\tilde v}^{k-3}+...+u{\tilde v}^{k-2}w)]=$$
$$=\alpha_k (q)[{\tilde v}^{k} +a{\tilde v}^{k-1}-\frac{q^2+1}{q^2}
(1+q^{-2} +...+q^{-2(k-2)})uw{\tilde v}^{k-2}]=$$
$$=\alpha_k (q)[{\tilde v}^{k} +a{\tilde v}^{k-1}+\frac{q^{-2(k-1)} -1}
{q^2 (q^2+1)(q^{-2} -1)}({\tilde v}^2+a(q^2+1){\tilde v}-{\tilde c}q^2)
{\tilde v}^{k-2}]=$$
$$=\alpha_k (q)[{\tilde v}^{k} +a{\tilde v}^{k-1}+\frac{q^{2(k-1)} -1}
{q^{2(k-1)}(q^2+1)(q^{2} -1)}({\tilde v}^k+a(q^2+1){\tilde v}^{k-1}-
{\tilde c}q^2 {\tilde v}^{k-2})]=$$
$$=\beta_k (q)[\frac{q^{2(k+1)} -1}{q^{2} -1}{\tilde v}^k +a(q^2+1)
\frac{q^{2k} -1}{q^{2} -1}{\tilde v}^{k-1} -{\tilde c}q^2
\frac{q^{2(k-1)} -1}{q^{2} -1}{\tilde v}^{k-2}],$$
with 
$$\beta_k (q)=\frac{\alpha_k (q)}{q^{2(k-1)}(q^{2} +1)}=
\frac{q^{2k}-1}{q^{2k-1}(q^{2} -1)}$$
(we assume that ${\tilde v}^{-1}={\tilde v}^{-2}=0$).

Thus, we have established the following
\begin{proposition} 
$$K{\tilde v}^{k}=q\beta_k (q)
[\frac{q^{2(k+1)} -1}{q^{2} -1}{\tilde v}^k +a(q^2+1)
\frac{q^{2k} -1}{q^{2} -1}{\tilde v}^{k-1} -{\tilde c}q^2
\frac{q^{2(k-1)} -1}{q^{2} -1}{\tilde v}^{k-2}].$$
\end{proposition}

\begin{remark}
This Proposition generalizes Proposition 6.2 from \cite{VS}. However, in order to 
represent it in a form similar to that from \cite{VS}, let us introduce the notions 
of right and left {\it q-difference} for a function $f(z),\ (z\in {\bf C})$ as follows:
$${\delta^+}_q f(z):=\frac{f(z)-f(qz)}{z-qz},\ \
{\delta^-}_q f(z):=\frac{f(z)-f(q^{-1}z)}{z-q^{-1}z}.$$

In particular, for $f(z)=z^k$ we have:
$${\delta^+}_{q^2} z^k = z^{k-1}\frac{q^{2k}-1}{q^2-1},\ \
{\delta^-}_{q^2} z^k = z^{k-1}\frac{q^{2k}-1}{q^{2(k-1)}(q^2-1)}.$$
Using this notation, we can rewrite the above formula for the action of the Casimir as follows:
$$K{\tilde v}^k =\frac{q^{2k} -1}{q^{2(k-1)}(q^{2} -1)}{\delta^+}_{q^2}
[({\tilde v}^2+a(q^2+1){\tilde v}-{\tilde c}q^2){\tilde v}^{k-1}]=$$
$$=[{\delta^+}_{q^2}({\tilde v}^2+a(q^2+1){\tilde v}-{\tilde c}q^2)
 {\delta^-}_{q^2}]{\tilde v}^k=
{\delta^+}_{q^2}({\tilde v}^2+h\frac{1+q^2}{1-q^2}{\tilde v}-
{\tilde c}q^2){\delta^-}_{q^2}{\tilde v}^k.$$

Thus, the action of the Casimir operator on any polynomial in ${\tilde v}$ can 
be expressed in terms of the $q^2$-difference operator of the second order.
\end{remark}

\section{$(\hbar,q)$-special polynomials}

It is well known that the function algebra $\aaa$ on a usual hyperboloid 
considered as an algebraic variety in $sl(2)^*$ is a direct sum of all integer 
spin irreducible $sl(2)$-modules ${\bf V}_k$. This property is also valid for 
its non-commutative analogue $\ah$.  It is also true if 
$q$ is generic for a $\uqs$-invariant algebra $\ahq$. 

To show this it suffices to check that for any integer spin $k$ there exists 
in the algebra $\ahq$ a unique polynomial in  ${\tilde v}$ belonging to the module 
${\bf V}_k$ (cf. \cite{DG}, where another method of the proof is given). 

This property is ensured by the following
\begin{proposition}
For any $\la_k,\, k=0,1,2,...$ given by the formulae (1) and for a generic $q$
there exists a unique polynomial of the form
$$P_k({\tilde v})=\sum^{k}_{j=0} A^{k}_j{\tilde v}^{k-j}\,\,\,
with\,\,\,A^{k}_0=1,$$
such that
$$K\, P_k({\tilde v})=\la_k\, P_k({\tilde v}).$$
\end{proposition}

Proof. Let $P_k({\tilde v})$ be such a polynomial. 

Using Proposition 1 we have:
$$KP_k({\tilde v})=\sum^{k}_{j=0} A^{k}_j(a_{k-j}{\tilde v}^{k-j}+b_{k-j}
{\tilde v}^{k-j-1}+c_{k-j}{\tilde v}^{k-j-2})=$$
$$a_k{\tilde v}^k+(b_k+A^{k}_1 a_{k-1}){\tilde v}^{k-1}+
\sum^{k}_{j=2}(A^{k}_{j-2}c_{k-j+2}+A^{k}_{j-1} b_{k-j+1}+A^{k}_j a_{k-j})
{\tilde v}^{k-j},$$
where 
$$a_k=\frac{(q^{2k}-1)(q^{2(k+1)}-1)}{q^{2k-2}(q^2-1)^2},\,
b_k=\frac{a(q^{2}+1)(q^{2k}-1)^2}{{q^{2k-2}(q^2-1)^2}},\,
c_k=-\frac{{\tilde c}(q^{2k}-1)(q^{2(k-1)}-1)}{q^{2k-4}(q^2-1)^2}.$$

Now the above equality of polynomials gives us the following chain of relations:
$$a_k=\lambda_k,$$
$$b_k+A^{k}_1 a_{k-1}=A^{k}_1 \lambda_{k},$$
$$A^{k}_{j-2}c_{k-j+2}+A^{k}_{j-1} b_{k-j+1}+A^{k}_j a_{k-j}=A^{k}_j
\lambda_{k} \,\, (j=2,3,...,k).$$

So we have the following recurrence relations for finding
$A^{k}_j$:
$$A^{k}_1=\frac{b_k}{a_k-a_{k-1}},$$
$$A^{k}_j=\frac{A^{k}_{j-2}c_{k-j+2}+A^{k}_{j-1} b_{k-j+1}}{a_k-a_{k-j}},
 \,\, (j=2,3,...,k). $$

It remains to say that the numerators of these formulae are not equal to 0 for 
a generic $q$. This completes the proof.

Let us remark that this approach gives us a description of "non-generic" values of 
$q$: they are exactly such that the numerators of the 
above formulae vanish. It should be noted that these numerators do not contain
$\hbar$ and therefore the decomposition $\ahq=\oplus {\bf V}_k$ is valid
for a generic $q$ independently on $\hbar$.

We call the above polynomials $(\hbar, q)$-{\em special polynomials}. If 
$\hbar=0$ and $q=1$, they coincide with the Legendre polynomials up to a change
of the variable and up to factors. A change of the variable consisting in 
multiplying the variable by $\sqrt{-1}$ is motivated by the fact that the 
Legendre polynomials arise from the real compact form of the group $SL(2,{\bf C})$. 

Since the Legendre polynomials are even for even $k$ and odd for odd $k$, 
this substitution leads to polynomials with real coefficients (for an odd $k$ 
it is necesary also to multiply the polynomial by $\sqrt{-1}$). 

It is still true if $\hbar=0$ but $q\not=1$. Thus, assuming $q$ to be real, in a similar way 
we get the polynomials with real coefficients which differ from the $q$-Legendre polynomials by 
factors (see \cite{VS},\cite{KN},\cite{Va}) .

However, if $\hbar\not=0$ and $q\not=1$, the above property is no longer true 
and the mentioned procedure does not lead to polynomials with real coefficients. 

This is the reason why we do prefer to deal with the complex form of the
quantum hyperboloid.  

\section{$(\hbar, q)$-integral}

Let us introduce in the algebra 
$\ahq$ an analogue of the invariant integral. It is exactly the projector in this algebra onto its
trivial component. 

In what follows we use the notation $\Int:\ahq \to {\bf C}$ for it. If $q=1,\,\hbar=0$, this 
operator coincides up to a factor with the usual invariant integral on a sphere.

In general case, we call this projector an $(\hbar, q)$-{\em integral}.

Our immediate aim is to compute the values Int$({\tilde v}^k)$. We use the
method of \cite{NM}, where a particular case  $(\hbar=0)$ has been considered.

It is obvious that $\Int (Yf)=0$ for any $f\in \ahq$.  This follows from the 
fact that $Yf\in{\bf V}_k$ if $f\in{\bf V}_k, k\neq 0$ and $Yf=0$ if $f\in {\bf V}_0$. 

Let us set $f=u{\tilde v}^k$. Then, using again the formula (2), we have:
$$\Int (Yu{\tilde
v}^k)=(YE_+^k+E_-YE_+^{k-1}+...+E_-^kY)u{\tilde v}^k=$$ 
$$-({\tilde v}+a){\tilde
v}^k+q^{-1}(q+q^{-1})u(w{\tilde v}^{k-1}+ {\tilde v}w {\tilde
v}^{k-2}+...+{\tilde v}^{k-1}w)=$$ 
$$-{\tilde v}^{k+1}-a{\tilde v}^k+(1+q^{-2})uw(1+q^{-2}+...+q^{-2(k-1)})
{\tilde v}^{k-1}=$$
$$-{\tilde v}^{k+1}-a{\tilde
v}^{k}+(q^2-1)^{-1}(1-q^{-2k})(1+q^2)^{-1} (-{\tilde v}^2-a(q^2+1){\tilde
v}+{\tilde c}q^2) {\tilde v}^{k-1}=0.$$

This implies the following equation
$$\mu_{k+1}(q^{2k+4}-1)+\mu_{k}a(q^{2k+2}-1)(1+q^2)-\mu_{k-1}(q^{2k}-1)q^2{\tilde
c}=0,$$
where $\mu_k=\Int({\tilde v}^{k})$. Now by putting
$\gamma_k=\mu_k(q^{2k+2}-1)$ we have
$$\gamma_{k+1}+a(1+q^2)\gamma_{k}-q^2{\tilde c}\gamma_{k-1}=0.$$

Thus, if we normalize the $(\hbar,q)$-integral by Int$(1)=1$ and Int$({\tilde
v})=0$, we have for $\mu_k$ the following formula:
\begin{equation}
\mu_k=(q^2-1)(q^{2k+2}-1)^{-1}(x_2x_1^k-x_1x_2^k)(x_2-x_1)^{-1},
\end{equation}
where $x_1$ and $x_2$ are the roots of the quadratic equation
$$x^2+a(1+q^2)x-q^2{\tilde c}^2=0.$$

Thus, we have proved the following
\begin{proposition} The $(\hbar,q)$-integral normalized by $\Int(1)=1$ and 
$\Int ({\tilde v})=0$ is unique and defined by the formula $\Int ({\tilde v}^k)=
\mu_k$, where $\mu_k$ is given by (3).

Assuming $\vert q\vert$ to be smaller than 1, we can represent this formula as:
$$\Int (f)=(1-q^2)(x_2-x_1)^{-1}\sum_{m=0}^{\infty}
(x_2f(x_1q^{2m})-x_1f(x_2q^{2m}))q^{2m},$$
where $f$ is a polynomial in ${\tilde v}$.
\end{proposition}

\begin{remark} a) Let us remark that $Int(P_k({\tilde v}))=0,\,k\geq 1$ for all
$(\hbar,q)$-special polynomials introduced in Section 4. Also, we have:
$$Int(P_k({\tilde v})P_l({\tilde v}))=0\ (k\not=l),$$  
i.e., $(\hbar,q)$-special polynomials are mutually orthogonal with respect 
to the pairing defined by the $(\hbar,q)$-integral. This follows from the fact that 
in the decomposition ${\bf V}_k\ot {\bf V}_l$ into a direct sum of irreducible 
components, the trivial component is present if and only if $k=l$.

b) Let us note that neither our formula for $(\hbar,q)$-special 
polynomials nor that for the $(\hbar,q)$-integral have any limit as $q\to 1$ if 
$\hbar\not=0$. In the classical case $(q=1)$ one usually deals with a family 
of finite dimensional representations of the algebras $\ah$. In such a 
representation, "the $(\hbar,1)$-integral" becomes a usual trace (up to a factor).

c) If $\hbar=0, q\not=0,$ the above formula for the $(\hbar,q)$-integral gives the well
known formula for the Jackson integral (see \cite{VS}, \cite{KN}). The relations of 
orthogonality for $(\hbar,q)$-special polynomials with respect to the 
$(\hbar,q)$-integral generalize those for $q$-Legendre polynomials with respect 
to the Jackson integral.
\end{remark}

\section{Connection with the deformation quantization}

By the fact that the algebras $\aaa$ and $\ahq$ are isomorphic as ${\bf C}[[\hbar,q]]$-modules 
(in other words, the deformation $\aaa\to\ahq$ is flat, here the parameters are assumed to be 
formal) one can introduce the corresponding quasiclassical object. This is {\em a Poisson pencil} 
(i.e., a linear space of Poisson brackets) generated by the KKS bracket and a so-called 
{\em R-matrix bracket} well defined on a hyperboloid.

The latter bracket is introduced by $\{f,g\}=<\rho^{\ot 2}(R), df\ot dg>$, where $R$ is 
the unique (up to a factor and an intertwinning) solution of the classical {\em modified} 
Yang-Baxter equation on the Lie algebra $sl(2)$, $\rho$ is the coadjoint representation 
restricted to a hyperboloid, and we use the pairing between vector fields and differential 
forms (extended onto their tensor powers)\footnote{As for other simple Lie algebras $g$ such 
a type of Poisson pencils exists only on some exeptional orbits in $g^*$ (cf. \cite{GP}).}.

Thus, the algebra $\ahq$ can be treated as a quantum object with respect to the above 
Poisson pencil. Let us emphasize that the quantization of the only KKS bracket gives the 
algebra $\ah$ which is $sl(2)$-invariant. Let us introduce an $sl(2)$-morphism 
$\phi:\aaa\to\ah$ by sending $u^k\in\aaa$ to $u^k\in\ah$.

By means of this morphism we can, in the spirit of the deformation quatization theory, 
introduce a new $sl(2)$-invariant associative product in the algebra $\aaa$: 
$$a\circ_{\hbar} b=\phi^{-1}(\phi(a)\circ\phi(b)),\,\, a,b\in\aaa,$$
where $\circ$ is the product in the algebra $\ah$. One can see that this quantization is 
closed in the sense of \cite{CFS} (this means that the trace in the quantum algebra is exactly 
an integral on the initial manifold, in our case such a manifold is a sphere). 

Let us remark that such a quantization exists for any symplectic Poisson bracket on any 
(compact smooth) manifold (cf. \cite{CFS}).

The passage $\aaa\to\ah$ is a particular case of this deformation quantization scheme since 
the KKS bracket is symplectic. It is not the case of the R-matrix bracket: it is not symplectic 
and its quantization leads to a deformation of the integral. Although it is easy, by  a 
method similar to the above, to represent the algebra $\ahq$ as $\aaa$ equipped with a 
deformed product $\circ_{\hbar,q}$, the initial integral on $\aaa$ is not any more a trace for 
this product.

Thus, the first step of the quantization, i.e., the passage $\aaa\to\ah$, can be done without 
any deformation of the integral. On the contrary, the second one, i.e., the further passage to 
the algebra $\ahq$, leads to such a deformation.

Let us explain now in what sense we use the term {\em q-commutative} for the algebra $\aq$. In 
this algebra there exists an involutive $(\tS^2=id)$ operator $\tS:(\aq)^{\ot 2}\to (\aq)^{\ot 2}$
which plays the role of an ordinary flip in the algebra $\aq$. It can be derived from the 
Yang-Baxter operator $S$: it suffices to replace all eigenvalues of $S$ close to 1 (resp. -1) by 1 
(resp. -1) keeping all eigenspaces of $S$ (assuming that $\vert q-1\vert\ll 1$). 

Another description of the operator $\tS$ is given in \cite{DS}. Using the results of this paper,
one can see that in the algebra $\aq$ we have $a\circ b=\circ(\tS(a\ot b))$ for any two elements
$a,b\in\aq$. In this sense, we say that the algebra $\aq$ is {\em q-commutative}.

Thus, quantizing the only R-matrix bracket, we pass from a commutative algebra to a q-commutative 
one. Meanwhile, a simultaneous quantization of the considered Poisson pencil leads to the algebras 
which are $U_q(\gggg)$-invariant but are no longer q-commutative. This gives a simultaneous
deformation of the category (instead of $sl(2)$-invariant algebras we get $\uq$-invariant ones) 
and a passage from "commutative" objects to "noncommutative" ones in the new category. 

We consider the final algebra $\ahq$ as an object of the twisted Quantum Mechanics, which looks 
like similar objects of the super-Quantum Mechanics. For a more detailed discussion of this point 
of view, we refer the reader to \cite{DGR1} and \cite{DGR2}.

{\bf Acknowledgements} This paper was begun during the stay of one of the authors (L.V.) at the 
Laboratory of Mathemathics of Valenciennes University. He would like to thank the University for 
a warm hospitality.

\end{document}